# High-Dimensional Stochastic Design Optimization by Adaptive-Sparse Polynomial Dimensional Decomposition

Sharif Rahman, Xuchun Ren, and Vaibhav Yadav

**Abstract** This paper presents a novel adaptive-sparse polynomial dimensional decomposition (PDD) method for stochastic design optimization of complex systems. The method entails an adaptive-sparse PDD approximation of a high-dimensional stochastic response for statistical moment and reliability analyses; a novel integration of the adaptive-sparse PDD approximation and score functions for estimating the first-order design sensitivities of the statistical moments and failure probability; and standard gradient-based optimization algorithms. New analytical formulae are presented for the design sensitivities that are simultaneously determined along with the moments or the failure probability. Numerical results stemming from mathematical functions indicate that the new method provides more computationally efficient design solutions than the existing methods. Finally, stochastic shape optimization of a jet engine bracket with 79 variables was performed, demonstrating the power of the new method to tackle practical engineering problems.

## 1 Introduction

Uncertainty quantification of complex systems, whether natural or man-made, is an important ingredient in numerous fields of science and engineering. For practical applications, encountering hundreds of input variables or more is not uncommon,

---

S. Rahman
College of Engineering and Program of Applied Mathematical & Computational Sciences, The University of Iowa, Iowa City, IA 52242, USA, e-mail: sharif-rahman@uiowa.edu

X. Ren
Department of Mechanical Engineering, Georgia Southern University, Statesboro, GA 30458, USA, e-mail: xren@georgiasouthern.edu

V. Yadav
Department of Aerospace Engineering, San Diego State University, San Diego, CA 92182, USA, e-mail: vyadav@mail.sdsu.edu

where an output function of interest, often defined algorithmically via finite-element analysis (FEA), is all too often expensive to evaluate. Modern surrogate methods, comprising stochastic collocation [1], polynomial chaos expansion [13], and sparse-grid quadrature [3], are known to offer significant computational advantages over crude Monte Carlo simulation (MCS). However, for truly high-dimensional systems, they require astronomically large numbers of terms or coefficients, succumbing to the curse of dimensionality. Therefore, alternative computational methods capable of exploiting low effective dimensions of multivariate functions, such as the polynomial dimensional decomposition (PDD) methods [7, 9], including a recently developed adaptive-sparse PDD method, are desirable [15]. Although PDD and PCE contain the same measure-consistent orthogonal polynomials, a recent work reveals that the error committed by the PDD approximation cannot be worse than that perpetrated by the PCE approximation for identical expansion orders [9].

An important application of uncertainty quantification is stochastic design optimization, which can be grouped in two principal classes: (1) design optimization for robustness [10], which minimizes the propagation of input uncertainty to output responses of interest, leading to an insensitive design; and (2) design optimization for reliability [11], which concentrates on attaining an optimal design by ensuring sufficiently low risk of failure. Depending on the objective set forth by a designer, uncertainty can be effectively mitigated by either class of design optimization. Indeed, with new formulations and methods appearing almost every year, stochastic design optimization in conjunction with FEA are becoming increasingly relevant and perhaps necessary for realistic design of complex structures and systems.

This paper presents an adaptive-sparse PDD method for stochastic design optimization of complex systems. The method is based on (1) an adaptive-sparse PDD approximation of a high-dimensional stochastic response for statistical moment and reliability analyses; (2) a novel integration of the adaptive-sparse PDD approximation and score functions for calculating the first-order sensitivities of the statistical moments and failure probability with respect to the design variables; and (3) standard gradient-based optimization algorithms. Section 2 formally defines two general variants of stochastic design optimization, including their concomitant mathematical statements. Section 3 starts with a brief exposition of the adaptive-sparse PDD approximation, leading to statistical moment and reliability analyses. Exploiting score functions, the section explains how the effort required to perform stochastic analyses also delivers the design sensitivities, sustaining no additional cost. The section also describes a coupling between stochastic analyses and design sensitivity analysis, resulting in an efficient optimization algorithm for solving both variants of the design optimization problem. Section 4 presents two numerical examples, including solving a large-scale shape design optimization problem. Finally, the conclusions are drawn in Sect. 5.

## 2 Stochastic Design Optimization

Consider a measurable space $(\Omega_{\mathbf{d}}, \mathscr{F}_{\mathbf{d}})$, where $\Omega_{\mathbf{d}}$ is a sample space and $\mathscr{F}_{\mathbf{d}}$ is a $\sigma$-field on $\Omega_{\mathbf{d}}$. Defined over $(\Omega_{\mathbf{d}}, \mathscr{F}_{\mathbf{d}})$, let $\{P_{\mathbf{d}} : \mathscr{F}_{\mathbf{d}} \to [0,1]\}$ be a family of probability measures, where for $M \in \mathbb{N} := \{1,2,\cdots\}$ and $N \in \mathbb{N}$, $\mathbf{d} = (d_1,\cdots,d_M) \in \mathscr{D}$ is an $\mathbb{R}^M$-valued design vector with non-empty closed set $\mathscr{D} \subseteq \mathbb{R}^M$ and let $\mathbf{X} := (X_1,\cdots,X_N) : (\Omega_{\mathbf{d}}, \mathscr{F}_{\mathbf{d}}) \to (\mathbb{R}^N, \mathscr{B}^N)$ be an $\mathbb{R}^N$-valued input random vector with $\mathscr{B}^N$ representing the Borel $\sigma$-field on $\mathbb{R}^N$, describing the statistical uncertainties in input variables of a complex system. The probability law of $\mathbf{X}$ is completely defined by a family of the joint probability density functions $\{f_{\mathbf{X}}(\mathbf{x};\mathbf{d}),\ \mathbf{x} \in \mathbb{R}^N,\ \mathbf{d} \in \mathscr{D}\}$ that are associated with probability measures $\{P_{\mathbf{d}},\ \mathbf{d} \in \mathscr{D}\}$, so that the probability triple $(\Omega_{\mathbf{d}}, \mathscr{F}_{\mathbf{d}}, P_{\mathbf{d}})$ of $\mathbf{X}$ depends on $\mathbf{d}$. A design variable $d_k$ can be any distribution parameter or a statistic — for instance, the mean or standard deviation — of $X_i$.

Let $y_l(\mathbf{X})$, $l = 0,1,\cdots,K$, be a collection of $K+1$ real-valued, square-integrable, measurable transformations on $(\Omega_{\mathbf{d}}, \mathscr{F}_{\mathbf{d}})$, describing performance functions of a complex system. It is assumed that $y_l : (\mathbb{R}^N, \mathscr{B}^N) \to (\mathbb{R}, \mathscr{B})$ is not an explicit function of $\mathbf{d}$, although $y_l$ implicitly depends on $\mathbf{d}$ via the probability law of $\mathbf{X}$. This is not a major limitation, as most design optimization problems involve means and/or standard deviations of random variables as design variables. There exist two prominent variants of design optimization under uncertainty: (1) design optimization for robustness and (2) design optimization for reliability. Their mathematical formulations, comprising an objective function $c_0 : \mathbb{R}^M \to \mathbb{R}$ and constraint functions $c_l : \mathbb{R}^M \to \mathbb{R}$, $l = 1,\cdots,K$, $1 \le K < \infty$, entail finding an optimal design solution $\mathbf{d}^*$ as follows.

- **Design for Robustness** [10]

$$
\begin{aligned}
\mathbf{d}^* = \underset{\mathbf{d} \in \mathscr{D} \subseteq \mathbb{R}^M}{\arg\min}\ c_0(\mathbf{d}) &:= w_1 \frac{\mathbb{E}_{\mathbf{d}}[y_0(\mathbf{X})]}{\mu_0^*} + w_2 \frac{\sqrt{\mathrm{var}_{\mathbf{d}}[y_0(\mathbf{X})]}}{\sigma_0^*}, \\
\text{subject to } c_l(\mathbf{d}) &:= \alpha_l \sqrt{\mathrm{var}_{\mathbf{d}}[y_l(\mathbf{X})]} - \mathbb{E}_{\mathbf{d}}[y_l(\mathbf{X})] \le 0,\ l = 1,\cdots,K,
\end{aligned}
\tag{1}
$$

where $\mathbb{E}_{\mathbf{d}}[y_l(\mathbf{X})] := \int_{\mathbb{R}^N} y_l(\mathbf{x}) f_{\mathbf{X}}(\mathbf{x};\mathbf{d}) d\mathbf{x}$ is the mean of $y_l(\mathbf{X})$ with $\mathbb{E}_{\mathbf{d}}$ denoting the expectation operator with respect to the probability measure $P_{\mathbf{d}}, \mathbf{d} \in \mathscr{D}$, $\mathrm{var}_{\mathbf{d}}[y_l(\mathbf{X})] := \mathbb{E}_{\mathbf{d}}[\{y_l(\mathbf{X}) - \mathbb{E}_{\mathbf{d}}[y_l(\mathbf{X})]\}^2]$ is the variance of $y_l(\mathbf{X})$, $w_1 \in \mathbb{R}_0^+ := [0,\infty)$ and $w_2 \in \mathbb{R}_0^+$ are two non-negative, real-valued weights, satisfying $w_1 + w_2 = 1$, $\mu_0^* \in \mathbb{R} \setminus \{0\}$ and $\sigma_0^* \in \mathbb{R}_0^+ \setminus \{0\}$ are two non-zero, real-valued scaling factors, and $\alpha_l \in \mathbb{R}_0^+$, $l = 0,1,\cdots,K$, are non-negative, real-valued constants associated with the probabilities of constraint satisfaction. For most applications, equal weights are chosen, but they can be distinct and biased, depending on the objective set forth by a designer. By contrast, the scaling factors are relatively arbitrary and chosen to better condition, such as normalize, the objective function. In (1), $c_0(\mathbf{d})$ describes the objective robustness, whereas $c_l(\mathbf{d})$, $l = 1,\cdots,K$, describe the feasibility robustness of a given design. The evaluations of both robustness measures involve the first two moments of various stochastic responses, consequently demanding statistical moment analysis.

- **Design for Reliability** [11]

$$
\begin{aligned}
\mathbf{d}^* &= \underset{\mathbf{d}\in\mathscr{D}\subseteq\mathbb{R}^M}{\arg\min}\; c_0(\mathbf{d}), \\
\text{subject to } c_l(\mathbf{d}) &:= P_{\mathbf{d}}\left[\mathbf{X} \in \Omega_{F,l}(\mathbf{d})\right] - p_l \le 0,\; l = 1,\cdots,K,
\end{aligned}
\tag{2}
$$

where $\Omega_{F,l}$ is the $l$th failure domain, $0 \le p_l \le 1$ is the $l$th target failure probability. In (2), the objective function $c_0$ is commonly prescribed as a deterministic function of $\mathbf{d}$, describing relevant system geometry, such as area, volume, and mass. In contrast, the constraint functions $c_l$, $l = 1,\cdots,K$, depending on the failure domain $\Omega_{F,l}$, require component or system reliability analyses. For a component reliability analysis, the failure domain is often adequately described by a single performance function $y_l(\mathbf{X})$, for instance, $\Omega_{F,l} := \{\mathbf{x} : y_l(\mathbf{x}) < 0\}$, whereas multiple, interdependent performance functions $y_{l,i}(\mathbf{x}),\; i = 1,2,\cdots$, are required for a system reliability analysis, leading, for example, to $\Omega_{F,l} := \{\mathbf{x} : \cup_i y_{l,i}(\mathbf{x}) < 0\}$ and $\Omega_{F,l} := \left\{\mathbf{x} : \cap_i y_{l,i}(\mathbf{x}) < 0\right\}$ for series and parallel systems, respectively.

The solution of a stochastic design optimization problem, whether in conjunction with robustness or reliability, mandates not only statistical moment and reliability analyses, but also the evaluations of gradients of moments and failure probability with respect to the design variables. The focus of this work is to solve a general high-dimensional design optimization problem described by (1) or (2) for arbitrary square-integrable functions $y_l(\mathbf{X})$, $l = 1,2,\cdots,K$, and for an arbitrary probability density $f_{\mathbf{X}}(\mathbf{x};\mathbf{d})$ of $\mathbf{X}$, provided that a few regularity conditions are met.

## 3 Adaptive-Sparse Polynomial Dimensional Decomposition Method

Let $y(\mathbf{X}) := y(X_1,\cdots,X_N)$ represent any one of the random functions $y_l$, $l = 0,1,\cdots,K$, introduced in Sect. 2 and let $\mathscr{L}_2(\Omega_{\mathbf{d}},\mathscr{F}_{\mathbf{d}},P_{\mathbf{d}})$ represent a Hilbert space of square-integrable functions $y$ with respect to the probability measure $f_{\mathbf{X}}(\mathbf{x};\mathbf{d})d\mathbf{x}$ supported on $\mathbb{R}^N$. Assuming independent coordinates, the joint probability density function of $\mathbf{X}$ is expressed by the product, $f_{\mathbf{X}}(\mathbf{x};\mathbf{d}) = \prod_{i=1}^{i=N} f_{X_i}(x_i;\mathbf{d})$, of marginal probability density functions $f_{X_i} : \mathbb{R} \to \mathbb{R}_0^+$ of $X_i$, each defined on its probability triple $(\Omega_{i,\mathbf{d}},\mathscr{F}_{i,\mathbf{d}},P_{i,\mathbf{d}})$ with a bounded or an unbounded support on $\mathbb{R}$, $i = 1,\cdots,N$. Then, for a given subset $u \subseteq \{1,\cdots,N\}$, $f_{\mathbf{X}_u}(\mathbf{x}_u;\mathbf{d}) := \prod_{p=1}^{|u|} f_{X_{i_p}}(x_{i_p};\mathbf{d})$ defines the marginal density function of the subvector $\mathbf{X}_u = \{X_{i_1},\cdots,X_{i_{|u|}}\}^T$ of $\mathbf{X}$.

Let $\{\psi_{ij}(X_i;\mathbf{d});\; j = 0,1,\cdots\}$ be a set of univariate orthonormal polynomial basis functions in the Hilbert space $\mathscr{L}_2(\Omega_{i,\mathbf{d}},\mathscr{F}_{i,\mathbf{d}},P_{i,\mathbf{d}})$ that is consistent with the probability measure $P_{i,\mathbf{d}}$ of $X_i$ for a given design $\mathbf{d}$, where $i = 1,\cdots,N$. For a given $u = \{i_1,\cdots,i_{|u|}\} \subseteq \{1,\cdots,N\}$, $1 \le |u| \le N$, $1 \le i_1 < \cdots < i_{|u|} \le N$, denote by $(\times_{p=1}^{p=|u|}\Omega_{i_p,\mathbf{d}}, \times_{p=1}^{p=|u|}\mathscr{F}_{i_p,\mathbf{d}}, \times_{p=1}^{p=|u|}P_{i_p,\mathbf{d}})$ the product probability triple of the subvector $\mathbf{X}_u$. Since the probability density function of $\mathbf{X}_u$ is separable (independent), the product polynomial $\psi_{u\mathbf{j}_{|u|}}(\mathbf{X}_u;\mathbf{d}) := \prod_{p=1}^{|u|} \psi_{i_p j_p}(X_{i_p};\mathbf{d})$, where $\mathbf{j}_{|u|} = (j_1,\cdots,j_{|u|}) \in \mathbb{N}_0^{|u|}$, $\mathbb{N}_0 := \mathbb{N} \cup \{0\}$, is a $|u|$-dimensional multi-index, constitutes an orthonormal basis in $\mathscr{L}_2(\times_{p=1}^{p=|u|}\Omega_{i_p,\mathbf{d}}, \times_{p=1}^{p=|u|}\mathscr{F}_{i_p,\mathbf{d}}, \times_{p=1}^{p=|u|}P_{i_p,\mathbf{d}})$.

The PDD of a square-integrable function $y$ represents a hierarchical expansion [7, 9]

$$y(\mathbf{X}) = y_\emptyset(\mathbf{d}) + \sum_{\emptyset \neq u \subseteq \{1,\cdots,N\}} \sum_{\substack{\mathbf{j}_{|u|} \in \mathbb{N}_0^{|u|} \\ j_1,\cdots,j_{|u|} \neq 0}} C_{u\mathbf{j}_{|u|}}(\mathbf{d}) \psi_{u\mathbf{j}_{|u|}}(\mathbf{X}_u;\mathbf{d}) \tag{3}$$

in terms of random multivariate orthonormal polynomials, where

$$y_\emptyset(\mathbf{d}) = \int_{\mathbb{R}^N} y(\mathbf{x}) f_{\mathbf{X}}(\mathbf{x};\mathbf{d}) d\mathbf{x} \tag{4}$$

and

$$C_{u\mathbf{j}_{|u|}}(\mathbf{d}) := \int_{\mathbb{R}^N} y(\mathbf{x}) \psi_{u\mathbf{j}_{|u|}}(\mathbf{x}_u;\mathbf{d}) f_{\mathbf{X}}(\mathbf{x};\mathbf{d}) d\mathbf{x},\ \emptyset \neq u \subseteq \{1,\cdots,N\},\ \mathbf{j}_{|u|} \in \mathbb{N}_0^{|u|} \tag{5}$$

are various expansion coefficients. The condition $j_1,\cdots,j_{|u|} \neq 0$ used in (3) and equations throughout the remainder of this paper implies that $j_k \neq 0$ for all $k = 1,\cdots,|u|$. Derived from the ANOVA dimensional decomposition [4], (3) provides an exact representation because it includes all main and interactive effects of input variables. For instance, $|u| = 0$ corresponds to the constant component function $y_\emptyset$, representing the mean effect of $y$; $|u| = 1$ leads to the univariate component functions, describing the main effects of input variables, and $|u| = S$, $1 < S \leq N$, results in the $S$-variate component functions, facilitating the interaction among at most $S$ input variables $X_{i_1},\cdots,X_{i_S}$, $1 \leq i_1 < \cdots < i_S \leq N$. Further details of PDD are available elsewhere [7, 9].

Equation (3) contains an infinite number of coefficients, emanating from infinite numbers of orthonormal polynomials. In practice, the number of coefficients must be finite, say, by retaining finite-order polynomials and reduced-degree interaction among input variables. For instance, an $S$-variate, $m$th-order PDD approximation [7, 9]

$$\tilde{y}_{S,m}(\mathbf{X}) = y_\emptyset(\mathbf{d}) + \sum_{\substack{\emptyset \neq u \subseteq \{1,\cdots,N\} \\ 1 \leq |u| \leq S}} \sum_{\substack{\mathbf{j}_{|u|} \in \mathbb{N}_0^{|u|}, ||\mathbf{j}_{|u|}||_\infty \leq m \\ j_1,\cdots,j_{|u|} \neq 0}} C_{u\mathbf{j}_{|u|}}(\mathbf{d}) \psi_{u\mathbf{j}_{|u|}}(\mathbf{X}_u;\mathbf{d}) \tag{6}$$

is generated, where $||\mathbf{j}_{|u|}||_\infty := \max(j_1,\cdots,j_{|u|})$ defines the $\infty$-norm and the integers $0 \leq S \leq N$ and $1 \leq m < \infty$ represent the largest degree of interactions among input variables and the largest order of orthogonal polynomials retained in a concomitant truncation of the sum in (3). It is important to clarify that the right side of (6) contains sums of at most $S$-dimensional PDD component functions of $y$. Therefore, the term "$S$-variate" used for the PDD approximation should be interpreted in the context of including at most $S$-degree interaction of input variables, even though $\tilde{y}_{S,m}$ is strictly an $N$-variate function. When $S \to N$ and $m \to \infty$, $\tilde{y}_{S,m}$ converges to $y$ in the mean-square sense, generating a hierarchical and convergent sequence of approximations [7, 9].

### 3.1 Adaptive-Sparse PDD Approximation

For practical applications, the dimensional hierarchy or nonlinearity of a stochastic response, in general, is not known *a priori*. Therefore, indiscriminately assigning values of the truncation parameters $S$ and $m$ is not desirable. Nor is it possible to do so when a stochastic solution is obtained via complex numerical algorithms. In which case, one must perform these truncations automatically by progressively drawing in higher-variate or higher-order contributions as appropriate. Based on the authors' past experience, an $S$-variate PDD approximation, where $S \ll N$, is adequate, when solving real-world engineering problems, with the computational cost varying polynomially ($S$-order) with respect to the number of variables [7, 9]. As an example, consider the selection of $S = 2$ for solving a stochastic problem in 100 dimensions by a bivariate PDD approximation, comprising $100 \times 99/2 = 4950$ bivariate component functions. If all such component functions are included, then the computational effort for even a full bivariate PDD approximation may exceed the computational budget allocated to solving this problem. But many of these component functions contribute little to the probabilistic characteristics sought and can be safely ignored. Similar conditions may prevail for higher-variate component functions. Henceforth, define an $S$-variate, partially adaptive-sparse PDD approximation [15]

$$\bar{y}_S(\mathbf{X}) := y_\emptyset(\mathbf{d}) + \sum_{\substack{\emptyset \neq u \subseteq \{1,\cdots,N\} \\ 1 \leq |u| \leq S}} \sum_{m_u=1}^{\infty} \sum_{\substack{\|\mathbf{j}_{|u|}\|_\infty = m_u,\, j_1,\cdots,j_{|u|} \neq 0 \\ \tilde{G}_{u,m_u} > \varepsilon_1, \Delta\tilde{G}_{u,m_u} > \varepsilon_2}} C_{u\mathbf{j}_{|u|}}(\mathbf{d}) \psi_{u\mathbf{j}_{|u|}}(\mathbf{X}_u; \mathbf{d}) \tag{7}$$

of $y(\mathbf{X})$, where

$$\tilde{G}_{u,m_u} := \frac{1}{\sigma^2(\mathbf{d})} \sum_{\substack{\mathbf{j}_{|u|} \in \mathbb{N}_0^{|u|},\, \|\mathbf{j}_{|u|}\|_\infty \leq m_u \\ j_1,\cdots,j_{|u|} \neq 0}} C_{u\mathbf{j}_{|u|}}^2(\mathbf{d}),\ m_u \in \mathbb{N},\ 0 < \sigma^2(\mathbf{d}) < \infty,$$

defines the approximate $m_u$th-order approximation of the global sensitivity index of $y(\mathbf{X})$ for a subvector $\mathbf{X}_u$, $\emptyset \neq u \subseteq \{1,\cdots,N\}$, of input variables $\mathbf{X}$ and

$$\Delta\tilde{G}_{u,m_u} := \begin{cases} \infty & \text{if } m_u = 1 \text{ or } (m_u \geq 2, \tilde{G}_{u,m_u-1} = 0, \tilde{G}_{u,m_u} \neq 0), \\ 0 & \text{if } m_u \geq 2, \tilde{G}_{u,m_u-1} = 0, \tilde{G}_{u,m_u} = 0, \\ \dfrac{\tilde{G}_{u,m_u} - \tilde{G}_{u,m_u-1}}{\tilde{G}_{u,m_u-1}} & \text{if } m_u \geq 2, \tilde{G}_{u,m_u-1} \neq 0 \end{cases}$$

defines the relative change in the approximate global sensitivity index when the largest polynomial order increases from $m_u - 1$ to $m_u$. The non-trivial definition applies when $m_u \geq 2$ and $\tilde{G}_{u,m_u-1} \neq 0$. When $m_u = 1$ or $(m_u \geq 2, \tilde{G}_{u,m_u-1} = 0, \tilde{G}_{u,m_u} \neq 0)$, the infinite value of $\Delta\tilde{G}_{u,m_u}$ guarantees that the $m_u$th-order contribution of $y_u$ to $y$ is preserved in the adaptive-sparse approximation. When $m_u \geq 2$, $\tilde{G}_{u,m_u-1} = 0$, and $\tilde{G}_{u,m_u} = 0$, the *zero* value of $\Delta\tilde{G}_{u,m_u}$ implies that there is no contribution of the $m_u$th-order contribution of $y_u$ to $y$. Here,

$$\sigma^2(\mathbf{d}) = \sum_{\emptyset \neq u \subseteq \{1,\cdots,N\}} \sum_{\substack{\mathbf{j}_{|u|} \in \mathbb{N}_0^{|u|} \\ j_1,\cdots,j_{|u|} \neq 0}} C^2_{u\mathbf{j}_{|u|}}(\mathbf{d}) \tag{8}$$

is the variance of $y(\mathbf{X})$. Then the sensitivity indices $\tilde{G}_{u,m_u}$ and $\Delta\tilde{G}_{u,m_u}$ provide an effective means to truncate the PDD in (3) both adaptively and sparsely. Equation (7) is attained by subsuming at most $S$-variate component functions, but fulfilling two inclusion criteria: (1) $\tilde{G}_{u,m_u} > \varepsilon_1$ for $1 \leq |u| \leq S \leq N$, and (2) $\Delta\tilde{G}_{u,m_u} > \varepsilon_2$ for $1 \leq |u| \leq S \leq N$, where $\varepsilon_1 \geq 0$ and $\varepsilon_2 \geq 0$ are two user-defined tolerances. The resulting approximation is partially adaptive because the truncations are restricted to at most $S$-variate component functions of $y$. When $S = N$, (7) becomes the fully adaptive-sparse PDD approximation [15]. The algorithmic details of numerical implementation associated with either fully or partially adaptive-sparse PDD approximation are available elsewhere [15].

The determination of PDD expansion coefficients $y_\emptyset(\mathbf{d})$ and $C_{u\mathbf{j}_{|u|}}(\mathbf{d})$ is vitally important for statistical moment and reliability analyses, including design sensitivities. As defined in (4) and (5), the coefficients involve $N$-dimensional integrals over $\mathbb{R}^N$. For large $N$, a multivariate numerical integration employing an $N$-dimensional tensor product of a univariate quadrature formula is computationally prohibitive and is, therefore, ruled out. An attractive alternative approach entails dimension-reduction integration [14], where the $N$-variate function $y$ in (4) and (5) is replaced by an $R$-variate $(1 \leq R \leq N)$ referential dimension decomposition at a chosen reference point. For instance, given a reference point $\mathbf{c} = (c_1,\cdots,c_N) \in \mathbb{R}^N$, the expansion coefficients $C_{u\mathbf{j}_{|u|}}$ are approximated by [14]

$$C_{u\mathbf{j}_{|u|}}(\mathbf{d}) \cong \sum_{i=0}^{R} (-1)^i \binom{N-R+i-1}{i} \sum_{\substack{v \subseteq \{1,\cdots,N\} \\ |v|=R-i, u \subseteq v}} \int_{\mathbb{R}^{|v|}} y(\mathbf{x}_v, \mathbf{c}_{-v}) \psi_{u\mathbf{j}_{|u|}}(\mathbf{x}_u;\mathbf{d}) f_{\mathbf{X}_v}(\mathbf{x}_v;\mathbf{d}) d\mathbf{x}_v, \tag{9}$$

requiring evaluations of at most $R$-dimensional integrals. The estimation of $y_\emptyset(\mathbf{d})$ is similar. The reduced integration facilitates calculation of the coefficients approaching their exact values as $R \to N$, and is significantly more efficient than performing one $N$-dimensional integration, particularly when $R \ll N$. Hence, the computational effort is significantly lowered using the dimension-reduction integration. For instance, when $R = 1$ or 2, (9) involves one-, or at most, two-dimensional integrations, respectively. Nonetheless, numerical integrations are still required for performing various $|v|$-dimensional integrals over $\mathbb{R}^{|v|}$, where $0 \leq |v| \leq R$. When $R > 1$, the multivariate integrals involved can be subsequently approximated by a sparse-grid quadrature, such as the fully symmetric interpolatory rule [5], as implemented by Yadav and Rahman [15].

## 3.2 Stochastic Analysis

### 3.2.1 Statistical Moments

Applying the expectation operator on $\bar{y}_S(\mathbf{X})$ and recognizing the *zero*-mean and orthogonal properties of PDD component functions, the mean

$$\mathbb{E}_{\mathbf{d}}\left[\bar{y}_S(\mathbf{X})\right] = y_\emptyset(\mathbf{d}) \tag{10}$$

of the partially adaptive-sparse PDD approximation agrees with the exact mean $\mathbb{E}_{\mathbf{d}}\left[y(\mathbf{X})\right] = y_\emptyset$ for any $\varepsilon_1$, $\varepsilon_2$, and $S$ [15]. However, the variance, obtained by applying the expectation operator on $(\bar{y}_S(\mathbf{X}) - y_\emptyset)^2$, varies according to [15]

$$\bar{\sigma}_S^2(\mathbf{d}) := \mathbb{E}_{\mathbf{d}}\left[\left(\bar{y}_S(\mathbf{X}) - \mathbb{E}\left[\bar{y}_S(\mathbf{X})\right]\right)^2\right] = \sum_{\substack{\emptyset \neq u \subseteq \{1,\cdots,N\} \\ 1 \leq |u| \leq S}} \sum_{m_u=1}^{\infty} \sum_{\substack{\|\mathbf{j}_{|u|}\|_\infty = m_u,\, j_1,\cdots,j_{|u|} \neq 0 \\ \tilde{G}_{u,m_u} > \varepsilon_1,\, \Delta\tilde{G}_{u,m_u} > \varepsilon_2}} C_{u\mathbf{j}_{|u|}}^2(\mathbf{d}), \tag{11}$$

where the squares of the expansion coefficients are summed following the same two pruning criteria discussed in the preceding subsection. Equation (11) provides a closed-form expression of the approximate second-moment properties of any square-integrable function $y$ in terms of the PDD expansion coefficients. When $S = N$ and $\varepsilon_1 = \varepsilon_2 = 0$, the right side of (11) coincides with that of (8). In consequence, the variance from the partially adaptive-sparse PDD approximation $\bar{y}_S(\mathbf{X})$ converges to the exact variance of $y(\mathbf{X})$ as $S \to N$, $\varepsilon_1 \to 0$, and $\varepsilon_2 \to 0$.

### 3.2.2 Failure Probability

A fundamental problem in reliability analysis entails calculation of the failure probability

$$P_F(\mathbf{d}) := P_{\mathbf{d}}\left[\mathbf{X} \in \Omega_F\right] = \int_{\mathbb{R}^N} I_{\Omega_F}(\mathbf{x}) f_{\mathbf{X}}(\mathbf{x};\mathbf{d}) d\mathbf{x} =: \mathbb{E}_{\mathbf{d}}\left[I_{\Omega_F}(\mathbf{X})\right],$$

where $I_{\Omega_F}(\mathbf{x})$ is the indicator function associated with the failure domain $\Omega_F$, which is equal to *one* when $\mathbf{x} \in \Omega_F$ and *zero* otherwise. Depending on component or system reliability analysis, let $\bar{\Omega}_{F,S} := \{\mathbf{x} : \bar{y}_S(\mathbf{x}) < 0\}$ or $\bar{\Omega}_{F,S} := \{\mathbf{x} : \cup_i \bar{y}_{i,S}(\mathbf{x}) < 0\}$ or $\bar{\Omega}_{F,S} := \{\mathbf{x} : \cap_i \bar{y}_{i,S}(\mathbf{x}) < 0\}$ be an approximate failure set as a result of $S$-variate, adaptive-sparse PDD approximations $\bar{y}_S(\mathbf{X})$ of $y(\mathbf{X})$ or $\bar{y}_{i,S}(\mathbf{X})$ of $y_i(\mathbf{X})$. Then the adaptive-sparse PDD estimate of the failure probability $P_F(\mathbf{d})$ is

$$\bar{P}_{F,S}(\mathbf{d}) = \mathbb{E}_{\mathbf{d}}\left[I_{\bar{\Omega}_{F,S}}(\mathbf{X})\right] = \lim_{L\to\infty}\frac{1}{L}\sum_{l=1}^{L} I_{\bar{\Omega}_{F,S}}(\mathbf{x}^{(l)}), \tag{12}$$

where $L$ is the sample size, $\mathbf{x}^{(l)}$ is the $l$th realization of $\mathbf{X}$, and $I_{\bar{\Omega}_{F,S}}(\mathbf{x})$ is another indicator function, which is equal to *one* when $\mathbf{x} \in \bar{\Omega}_{F,S}$ and *zero* otherwise.

Note that the simulation of the PDD approximation in (12) should not be confused with crude MCS commonly used for producing benchmark results. The crude MCS, which requires numerical calculations of $y(\mathbf{x}^{(l)})$ or $y_i(\mathbf{x}^{(l)})$ for input samples $\mathbf{x}^{(l)}$, $l = 1,\cdots,L$, can be expensive or even prohibitive, particularly when the sample size $L$ needs to be very large for estimating small failure probabilities. In contrast, the MCS embedded in the adaptive-sparse PDD approximation requires evaluations of simple polynomial functions that describe $\bar{y}_S(\mathbf{x}^{(l)})$ or $\bar{y}_{i,S}(\mathbf{x}^{(l)})$. Therefore, a relatively large sample size can be accommodated in the adaptive-sparse PDD method even when $y$ or $y_i$ is expensive to evaluate.

## 3.3 Design Sensitivity Analysis

When solving design optimization problems employing gradient-based optimization algorithms, at least the first-order derivatives of the first two moments and failure probability with respect to each design variable are required. In this subsection, the adaptive-sparse PDD method for statistical moment and reliability analyses is expanded for design sensitivity analysis. For such sensitivity analysis, the following regularity conditions are assumed: (1) The domains of design variables $d_k \in \mathscr{D}_k \subset \mathbb{R}$, $k = 1,\cdots,M$, are open intervals of $\mathbb{R}$; (2) the probability density function $f_{\mathbf{X}}(\mathbf{x};\mathbf{d})$ of $\mathbf{X}$ is continuous. In addition, the partial derivative $\partial f_{\mathbf{X}}(\mathbf{x};\mathbf{d})/\partial d_k$, $k = 1,\cdots,M$, exists and is finite for all $\mathbf{x} \in \mathbb{R}^N$ and $d_k \in \mathscr{D}_k$. Furthermore, the statistical moments of $y$ and failure probability are differentiable functions of $\mathbf{d} \in \mathscr{D} \subseteq \mathbb{R}^M$; and (3) there exists a Lebesgue integrable dominating function $z(\mathbf{x})$ such that $|y^r(\mathbf{x})\partial f_{\mathbf{X}}(\mathbf{x};\mathbf{d})/\partial d_k| \le z(\mathbf{x})$ and $|I_{\Omega_F}(\mathbf{x})\partial f_{\mathbf{X}}(\mathbf{x};\mathbf{d})/\partial d_k| \le z(\mathbf{x})$, where $r = 1,2$, and $k = 1,\cdots,M$.

### 3.3.1 Score Function

Let

$$h(\mathbf{d}) := \mathbb{E}_{\mathbf{d}}[g(\mathbf{X})] := \int_{\mathbb{R}^N} g(\mathbf{x}) f_{\mathbf{X}}(\mathbf{x};\mathbf{d})d\mathbf{x} \tag{13}$$

be a generic probabilistic response, where $h(\mathbf{d})$ and $g(\mathbf{x})$ are either the $r$th-order raw moment $m^{(r)}(\mathbf{d}) := \mathbb{E}_{\mathbf{d}}[y_S^r(\mathbf{X})]$ $(r = 1,2)$ and $y^r(\mathbf{x})$ for statistical moment analysis or $P_F(\mathbf{d})$ and $I_{\Omega_F}(\mathbf{x})$ for reliability analysis. Suppose that the first-order derivative of $h(\mathbf{d})$ with respect to a design variable $d_k$, $1 \le k \le M$, is sought. Taking the partial derivative of $h(\mathbf{d})$ with respect to $d_k$ and then applying the Lebesgue dominated convergence theorem [2], which permits the differential and integral operators to be interchanged, yields the first-order sensitivity

$$\begin{aligned}\frac{\partial h(\mathbf{d})}{\partial d_k} &:= \frac{\partial \mathbb{E}_{\mathbf{d}}\left[g(\mathbf{X})\right]}{\partial d_k} \\ &= \frac{\partial}{\partial d_k}\int_{\mathbb{R}^N} g(\mathbf{x}) f_{\mathbf{X}}(\mathbf{x};\mathbf{d}) d\mathbf{x} \\ &= \int_{\mathbb{R}^N} g(\mathbf{x}) \frac{\partial \ln f_{\mathbf{X}}(\mathbf{x};\mathbf{d})}{\partial d_k} f_{\mathbf{X}}(\mathbf{x};\mathbf{d}) d\mathbf{x} \\ &=: \mathbb{E}_{\mathbf{d}}\left[g(\mathbf{X}) s_{d_k}^{(1)}(\mathbf{X};\mathbf{d})\right],\end{aligned} \tag{14}$$

provided that $f_{\mathbf{X}}(\mathbf{x};\mathbf{d}) > 0$ and the derivative $\partial \ln f_{\mathbf{X}}(\mathbf{x};\mathbf{d})/\partial d_k$ exists. In the last line of (14), $s_{d_k}^{(1)}(\mathbf{X};\mathbf{d}) := \partial \ln f_{\mathbf{X}}(\mathbf{X};\mathbf{d})/\partial d_k$ is known as the first-order score function for the design variable $d_k$ [8, 12]. According to (13) and (14), the generic probabilistic response and its sensitivities have both been formulated as expectations of stochastic quantities with respect to the same probability measure, facilitating their concurrent evaluations in a single stochastic simulation or analysis.

### 3.3.2 Sensitivity of Statistical Moments

Selecting $h(\mathbf{d})$ and $g(\mathbf{x})$ to be $m^{(r)}(\mathbf{d})$ and $y^r(\mathbf{x})$, respectively, and then replacing $y(\mathbf{x})$ with its $S$-variate adaptive-sparse PDD approximation $\bar{y}_S(\mathbf{x})$ in the last line of (14), the resultant approximation of the sensitivities of the $r$th-order moment is obtained as

$$\mathbb{E}_{\mathbf{d}}\left[\bar{y}_S^r(\mathbf{X}) s_{d_k}^{(1)}(\mathbf{X};\mathbf{d})\right] = \int_{\mathbb{R}^N} \bar{y}_S^r(\mathbf{x}) s_{d_k}^{(1)}(\mathbf{x};\mathbf{d}) f_{\mathbf{X}}(\mathbf{x};\mathbf{d}) d\mathbf{x}. \tag{15}$$

The $N$-dimensional integral in (15) can be estimated by the same or similar dimension-reduction integration as employed for estimating the PDD expansion coefficients. Furthermore, if $s_{d_k}^{(1)}$ is square-integrable, then it can be expanded with respect to the same orthogonal polynomial basis functions, resulting in a closed-form expression of the design sensitivity [8]. Finally, setting $r = 1$ or 2 in (15) delivers the approximate sensitivity of the first or second moment.

### 3.3.3 Sensitivity of Failure Probability

Selecting $h(\mathbf{d})$ and $g(\mathbf{x})$ to be $P_F(\mathbf{d})$ and $I_{\Omega_F}(\mathbf{x})$, respectively, and then replacing $y(\mathbf{x})$ with its $S$-variate adaptive-sparse PDD approximation $\bar{y}_S(\mathbf{x})$ in the last line of (14), the resultant approximation of the sensitivities of the failure probability is obtained as

$$\mathbb{E}_{\mathbf{d}}\left[I_{\bar{\Omega}_{F,S}}(\mathbf{X}) s_{d_k}^{(1)}(\mathbf{X};\mathbf{d})\right] = \lim_{L\to\infty}\frac{1}{L}\sum_{l=1}^{L}\left[I_{\bar{\Omega}_{F,S}}(\mathbf{x}^{(l)}) s_{d_k}^{(1)}(\mathbf{x}^{(l)};\mathbf{d})\right], \tag{16}$$

where $L$ is the sample size, $\mathbf{x}^{(l)}$ is the $l$th realization of $\mathbf{X}$, and $I_{\bar{\Omega}_{F,S}}(\mathbf{x})$ is the adaptive-sparse PDD-generated indicator function. Again, the sensitivity in (16) is easily and inexpensively determined by sampling elementary polynomial functions that describe $\bar{y}_S$ and $s_{d_k}^{(1)}$.

*Remark 1.* The PDD expansion coefficients depend on the design vector $\mathbf{d}$. Naturally, a PDD approximation, whether obtained by truncating arbitrarily or adaptively, is also dependent on $\mathbf{d}$, unless the approximation exactly reproduces the function $y(\mathbf{X})$. It is important to clarify that the approximate sensitivities in (15) and (16) are obtained not by taking partial derivatives of the approximate moments in (10) and (11) and approximate failure probability in (12) with respect to $d_k$. Instead, they result from replacing $y(\mathbf{x})$ with $\bar{y}_S(\mathbf{x})$ in the expectation describing the last line of (14).

*Remark 2.* The score function method has the nice property that it requires differentiating only the underlying probability density function $f_{\mathbf{X}}(\mathbf{x};\mathbf{d})$. The resulting score functions can be easily and, in most cases, analytically determined. If the performance function is not differentiable or discontinuous − for example, the indicator function that comes from reliability analysis − the proposed method still allows evaluation of the sensitivity if the density function is differentiable. In reality, the density function is often smoother than the performance function, and therefore the proposed sensitivity methods are able to calculate sensitivities for a wide variety of complex mechanical systems.

### 3.4 Optimization Algorithm

The adaptive-sparse PDD approximations described in the preceding subsections provide a means to evaluate the objective and constraint functions, including their design sensitivities, from a single stochastic analysis. An integration of statistical moment analysis, reliability analysis, design sensitivity analysis, and a suitable optimization algorithm should render a convergent solution of the design optimization problems in (1) or (2). Algorithm 1 describes the computational flow of the adaptive-sparse PDD method for stochastic design optimization.

## 4 Numerical Examples

Two examples are presented to illustrate the adaptive-sparse PDD method for design optimization under uncertainty, where the objective and constraint functions are either elementary mathematical constructs or relate to complex engineering problems. Orthonormal polynomials consistent with the probability distributions of input random variables were used as bases. The PDD expansion coefficients were estimated using dimension-reduction integration and sparse-grid quadrature entailing an extended fully symmetric interpolatory rule [5, 15]. The sensitivities of moments and

**Algorithm 1:** Proposed adaptive-sparse PDD for stochastic design optimization

```
Input: an initial design d_0, S, ε > 0, ε_1 > 0, ε_2 > 0, q = 0
Output: an approximation d*_S of optimal design d*
d^(q) ← d_0
repeat
    d ← d^(q)
    Generate adaptive-sparse PDD approximations ȳ_{l,S}(X) at current design d of all y_l(X)
    in (1) or (2), l = 0, 1, ··· , K
    if design for robustness then
        compute moments E_d[ȳ_{l,S}(X)] and σ̄²_{l,S}(d) of ȳ_{l,S}(X)
                                                    /* from (10) and (11) */
        estimate design sensitivity of moments
                                                    /* from (15) */
    else if design for reliability then
        compute failure probability P̄_{F,S}(d) for ȳ_{l,S}(X)
                                                    /* from (12) */
        estimate design sensitivity of failure probability
                                                    /* from (16) */
    endif
    Evaluate objective and constraint functions in (1) or (2) and their sensitivities at d
    Using a gradient-based algorithm, obtain the next design d^(q+1)
    Set q = q + 1
until ||d^(q) − d^(q−1)||_2 < ε
d*_S ← d^(q)
```

failure probability were evaluated using dimension-reduction integration and embedded MCS of the adaptive-sparse PDD approximation, respectively. The optimization algorithm selected is sequential quadratic programming in both examples.

## *4.1 Example 1: Mathematical Functions*

The first example entails design optimization for robustness, which calls for finding

$$\begin{aligned} \mathbf{d}^* = \underset{\mathbf{d}\in\mathscr{D}\subseteq\mathbb{R}^M}{\arg\min}\; & c_0(\mathbf{d}) := 0.5\frac{\mathbb{E}_{\mathbf{d}}\left[y_0(\mathbf{X})\right]}{10} + 0.5\frac{\sqrt{\mathrm{var}_{\mathbf{d}}\left[y_0(\mathbf{X})\right]}}{2}, \\ \text{subject to}\; & c_1(\mathbf{d}) := 3\sqrt{\mathrm{var}_{\mathbf{d}}\left[y_1(\mathbf{X})\right]} - \mathbb{E}_{\mathbf{d}}\left[y_1(\mathbf{X})\right] \le 0, \\ & c_2(\mathbf{d}) := 3\sqrt{\mathrm{var}_{\mathbf{d}}\left[y_2(\mathbf{X})\right]} - \mathbb{E}_{\mathbf{d}}\left[y_2(\mathbf{X})\right] \le 0, \end{aligned}$$

where $\mathbf{d} \in \mathscr{D} = [0.00002, 0.002] \times [0.1, 1.6]$ and

$$y_0(\mathbf{X}) = X_3 X_1 \sqrt{1 + X_2^2}$$

and

$$y_l(\mathbf{X}) = 1 - \frac{5X_4\sqrt{1+X_2^2}}{\sqrt{65}X_5}\left(\frac{8}{X_1} + (-1)^{l+1}\frac{1}{X_1X_2}\right),\ l = 1,2,$$

are three random response functions of five independent random variables. The first two variables $X_1$ and $X_2$ follow Gaussian distributions with respective means $d_1 = \mathbb{E}_{\mathbf{d}}[X_1]$ and $d_2 = \mathbb{E}_{\mathbf{d}}[X_2]$ and coefficients of variations both equal to 0.02. The remaining variables, $X_3$, $X_4$, and $X_4$, follow Beta, Gumbel, and Lognormal distributions with respective means of 10,000, 0.8, and 1050 and respective coefficients of variations of 0.2, 0.25 and 0.238. The initial design vector is $\mathbf{d}_0 = (0.001, 1)$. In this example, $N = 5$, $M = 2$, and $K = 2$.

Table 1 presents detailed optimization results from two distinct adaptive-sparse PDD approximations, entailing univariate ($S = 1$) and bivariate ($S = 2$) truncations, employed to solve this optimization problem. The optimal solutions by these two approximations are close to each other, both indicating that the first constraint is nearly active ($c_1 \cong 0$). The results of the bivariate approximations confirm that the univariate solution is adequate. However, the total numbers of function evaluations step up for the bivariate approximation, as expected.

Since this problem can also be solved by the non-adaptively truncated PDD [10] and tensor product quadrature (TPQ) [6] methods, comparing their solutions, listed in the last three columns of Table 1, with the adaptive-sparse PDD solutions should be intriguing. It appears that the existing PDD truncated at the largest polynomial order of the adaptive-sparse PDD approximation, which is four, and TPQ methods are also capable of producing a similar optimal solution, but by incurring computational cost far more than the adaptive-sparse PDD methods. For instance, comparing the total numbers of function evaluations, the univariate adaptive-sparse PDD method is more economical than the existing univariate PDD and TPQ methods by factors of 1.5 and 37.7, respectively. The new bivariate adaptive-sparse PDD is more than twice as efficient as the existing non-adaptively truncated bivariate PDD.

**Table 1** Optimization results for Example 1

| | Adaptive-Sparse PDD | | Truncated PDD [10] | | |
|---|---|---|---|---|---|
| Results | Univariate | Bivariate | Univariate | Bivariate | TPQ [6] |
| $d_1^*$ ($\times 10^{-4}$) | 11.3902 | 11.5753 | 11.3921 | 11.5695 | 11.6476 |
| $d_2^*$ | 0.3822 | 0.3780 | 0.3817 | 0.3791 | 0.3767 |
| $c_0(\mathbf{d}^*)$ | 1.2226 | 1.2408 | 1.2227 | 1.2406 | 1.2480 |
| $c_1(\mathbf{d}^*)$ | 0.0234 | 0.0084 | 0.0233 | 0.0084 | 0.0025 |
| $c_2(\mathbf{d}^*)$ | -0.4810 | -0.4928 | -0.4816 | -0.4917 | -0.4970 |
| No. of iterations | 12 | 13 | 12 | 14 | 10 |
| Total no. of function evaluations | 465 | 2374 | 696 | 6062 | 17,521 |

### *4.2 Example 2: Shape Optimization of a Jet Engine Bracket*

The final example demonstrates the usefulness of the adaptive-sparse PDD method in designing for reliability an industrial-scale mechanical component, known as jet engine bracket, as shown in Fig. 1(a). Seventy-nine random shape parameters, $X_i$, $i = 1, \cdots, 79$, resulting from manufacturing variability, describe the shape of a jet engine bracket in three dimensions, including two quadrilateral holes introduced to reduce the mass of the jet engine bracket as much as possible. The design variables, $d_i = \mathbb{E}_{\mathbf{d}}[X_i]$, $i = 1, \cdots, 79$, are the means of these 79 independent random variables, with Figures 1(b)-(d) depicting the initial design of the jet engine bracket geometry at mean values of the shape parameters. The centers of the four bottom circular holes are fixed; a deterministic force, $F = 43.091$ kN, was applied at the center of the top circular hole with a 48° angle from the horizontal line, as shown in Fig. 1(c), and a deterministic torque, $T = 0.1152$ kN-m, was applied at the center of the top circular hole, as shown in Fig. 1(d). The jet engine bracket is made of Titanium Alloy Ti-6Al-4V with deterministic material properties described elsewhere [11]. Due to their finite bounds, the random variables $X_i$, $i = 1, \cdots, 79$, were assumed to follow truncated Gaussian distributions [11].

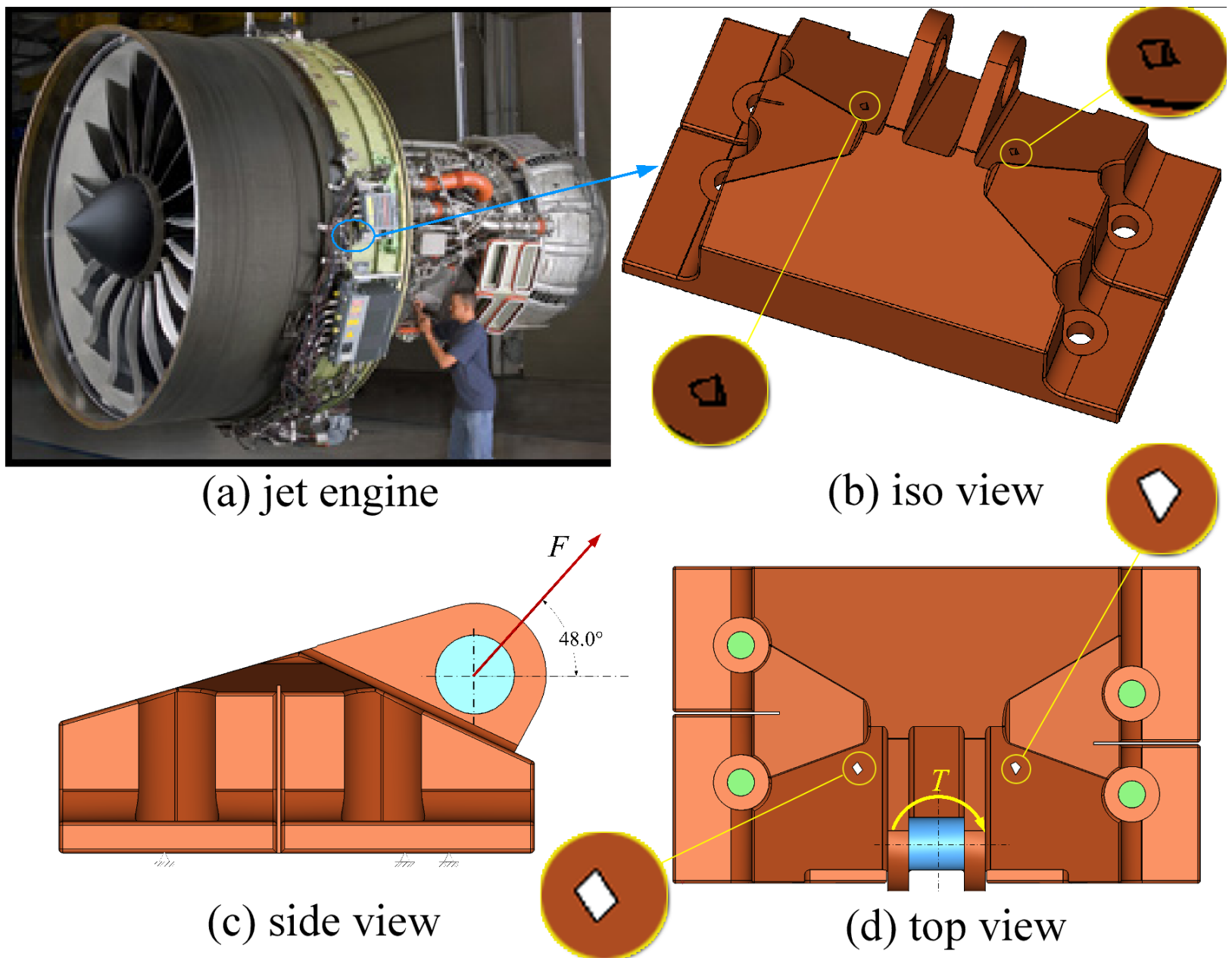

**Fig. 1** A jet engine bracket; (a) a jet engine; (b) isometric view; (c) lateral view; (d) top view.

The objective of this example is to minimize the mass of the engine bracket by changing the shape of the geometry such that the fatigue life $y(\mathbf{u}(\boldsymbol{\xi};\mathbf{X}), \boldsymbol{\sigma}(\boldsymbol{\xi};\mathbf{X}))$ ex-ceeds a million loading cycles with 99.865% probability. The underlying stochastic differential equations call for finding the displacement $\mathbf{u}(\boldsymbol{\xi};\mathbf{X})$ and stress $\boldsymbol{\sigma}(\boldsymbol{\xi};\mathbf{X})$ solutions at a spatial coordinate $\boldsymbol{\xi} = (\xi_1, \xi_2, \xi_3) \in \Omega \subset \mathbb{R}^3$, satisfying $P_{\mathbf{d}}$-almost surely

$$\begin{aligned} \nabla \cdot \boldsymbol{\sigma}(\boldsymbol{\xi};\mathbf{X}) + \mathbf{b}(\boldsymbol{\xi};\mathbf{X}) &= \mathbf{0} \text{ in } \Omega \subset \mathbb{R}^3, \\ \boldsymbol{\sigma}(\boldsymbol{\xi};\mathbf{X}) \cdot \mathbf{n}(\boldsymbol{\xi};\mathbf{X}) &= \bar{\mathbf{t}}(\boldsymbol{\xi};\mathbf{X}) \text{ on } \partial\Omega_t, \\ \mathbf{u}(\boldsymbol{\xi};\mathbf{X}) &= \bar{\mathbf{u}}(\boldsymbol{\xi};\mathbf{X}) \text{ on } \partial\Omega_u, \end{aligned} \tag{17}$$

such that $\partial\Omega_t \cup \partial\Omega_u = \partial\Omega$ and $\partial\Omega_t \cap \partial\Omega_u = \emptyset$ with $\nabla := (\partial/\partial\xi_1, \partial/\partial\xi_2, \partial/\partial\xi_3)$ and $\mathbf{b}(\boldsymbol{\xi};\mathbf{X})$, $\bar{\mathbf{t}}(\boldsymbol{\xi};\mathbf{X})$, $\bar{\mathbf{u}}(\boldsymbol{\xi};\mathbf{X})$, and $\mathbf{n}(\boldsymbol{\xi};\mathbf{X})$ representing the body force, prescribed traction on $\partial\Omega_t$, prescribed displacement on $\partial\Omega_u$, and unit outward normal vector, respectively. Mathematically, the problem entails finding an optimal design solution

$$\begin{aligned} \mathbf{d}^* = \underset{\mathbf{d}\in\mathscr{D}\subseteq\mathbb{R}^{79}}{\arg\min}\; c_0(\mathbf{d}) := \rho \int_{\Omega(\mathbf{d})} d\Omega \\ \text{subject to } c(\mathbf{d}) := P_{\mathbf{d}}\left[y_{\min}(\mathbf{u}(\boldsymbol{\xi}_c;\mathbf{X}), \boldsymbol{\sigma}(\boldsymbol{\xi}_c;\mathbf{X})) \le 10^6\right] \le 1 - 0.99865, \end{aligned}$$

where the objective function $c_0(\mathbf{d})$, with $\rho$ representing the mass density of the material, describes the overall mass of the bracket; on the other hand, the constraint function $c(\mathbf{d})$ quantifies the probability of minimum fatigue crack-initiation life $y_{\min}$, attained at a critical spatial point $\boldsymbol{\xi}_c$, failing to exceed a million loading cycles to be less than $(1 - 0.99865)$. Here, $y_{\min}$ depends on displacement and stress responses $\mathbf{u}(\boldsymbol{\xi}_c;\mathbf{X})$ and $\boldsymbol{\sigma}(\boldsymbol{\xi}_c;\mathbf{X})$, which satisfy (17). An FEA comprising 341,112 nodes and 212,716 ten-noded, quadratic, tetrahedral elements, was performed to solve the variational weak form of (17). Further details are available elsewhere [11].

The univariate ($S = 1$) adaptive-sparse PDD method was applied to solve this shape optimization problem. Figures 2(a) through (d) show the contour plots of the logarithm of fatigue life at mean shapes of several design iterations, including the initial design, throughout the optimization process. Due to a conservative initial design, with fatigue life contour depicted in Fig. 2(a), the minimum fatigue crack-initiation life of $6.65 \times 10^9$ cycles is much larger than the required fatigue crack-initiation life of a million cycles. For the tolerance and subregion size parameters selected, 14 iterations and 2,808 FEA led to a final optimal design with the corresponding mean shape presented in Fig. 2(d). The total run time, including performing all 2808 FEA in a desktop personal computer (8 cores, 2.3 GHz, 16 GB RAM), was about 165 hours. Most design variables have undergone significant changes from their initial values, prompting substantial modifications of the shapes or sizes of the outer boundaries, quadrilateral holes, and bottom surfaces of the engine bracket. The mean optimal mass of the engine bracket is 0.48 kg – an almost 84 percent reduction from the mean initial mass of 3.02 kg. At optimum, the constraint function $c(\mathbf{d})$ is practically zero and is, therefore, close to being active.

This example shows some promise of the adaptive-sparse PDD methods in solving industrial-scale engineering design problems with an affordable computational cost. However, an important drawback persists: given the computer resources available at the time of this work, only the univariate adaptive-sparse PDD approximation is feasible. The univariate result has yet to be verified with those obtained from bivariate or higher-variate adaptive-sparse PDD approximations. Therefore, the univariate “optimal” solution reported here should be guardedly interpreted.

Finally, it is natural to ask how much the bivariate adaptive-sparse PDD approximation will cost to solve this design problem. Due to quadratic computational complexity, the full bivariate PDD approximation using current computer resources of this study is prohibitive. However, a bivariate adaptive-sparse PDD approximation with a cost scaling markedly less than quadratic, if it can be developed, should be encouraging. In which case, a designer should exploit the univariate solution as the initial design to seek a better design using the bivariate adaptive-sparse PDD method, possibly, in fewer design iterations. The process can be repeated for higher-variate PDD methods if feasible. Clearly, additional research on stochastic design optimization, including more efficient implementation of the adaptive-sparse PDD methods, is required.

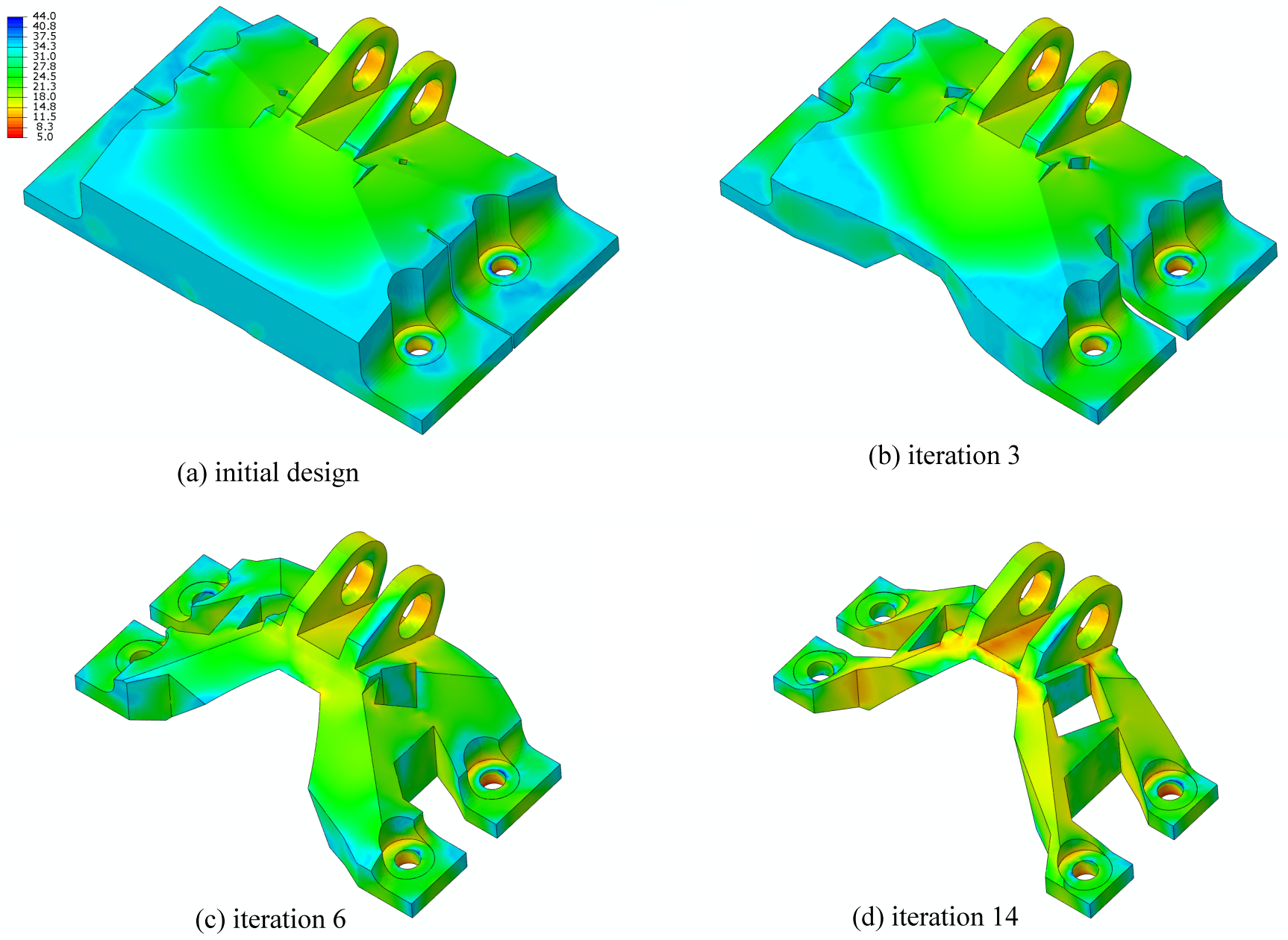

**Fig. 2** Contours of logarithmic fatigue life at mean shapes of the jet engine bracket by the adaptive-sparse PDD method; (a) initial design; (b) iteration 3; (c) iteration 6; (d) iteration 14 (optimum).

# 5 Conclusion

A new adaptive-sparse PDD method was developed for stochastic design optimization of high-dimensional complex systems commonly encountered in applied sciences and engineering. The method is based on an adaptive-sparse PDD approximation of a high-dimensional stochastic response for statistical moment and reliability analyses; a novel integration of the adaptive-sparse PDD approximation and score functions for estimating the first-order sensitivities of the statistical moments and failure probability with respect to the design variables; and standard gradient-based optimization algorithms, encompassing a computationally efficient design process. When blended with score functions, the adaptive-sparse PDD approximation leads to analytical formulae for calculating the design sensitivities. More importantly, the statistical moments, failure probability, and their respective design sensitivities are all determined concurrently from a single stochastic analysis or simulation. Numerical results stemming from a mathematical example indicate that the new method provides more computationally efficient design solutions than the existing methods. Finally, stochastic shape optimization of a jet engine bracket with 79 variables was performed, demonstrating the power of the new methods to tackle practical engineering problems.

**Acknowledgements** The authors acknowledge financial support from the U.S. National Science Foundation under Grant Nos. CMMI-0969044 and CMMI-1130147.